\documentclass[a4paper]{article}
\usepackage{etex}
\usepackage[latin1]{inputenc}
\usepackage{amsthm}
\usepackage{amssymb}
\usepackage{dsfont}         
\usepackage{graphicx}
\usepackage{tikz}
\usetikzlibrary{shapes,arrows}
\usepackage{epsf}           
\usepackage{times}

\newtheorem{theorem}{Theorem}[section]
\newtheorem{proposition}[theorem]{Proposition}

\newtheorem{theorem*}{Theorem}
\newtheorem{proposition*}{Proposition}
\newtheorem{corollary*}{Corollary}
\newtheorem{lemma*}{Lemma}

\newtheorem{example}[theorem]{Example}

\newtheorem{remark}[theorem]{Remark}

\newtheorem{definition*}{Definition}
\newtheorem{example*}{Example}
\newtheorem{examples*}{Examples}
\newtheorem{remark*}{Remark}
\newtheorem{remarks*}{Remarks}

\newcommand{\supp}{\mbox{\rm supp}}
\renewcommand{\epsilon}{\varepsilon}

\newcommand{\RR}{\mathbb{R}}

\newcommand{\NN}{\mathbb{N}}

\newcommand{\cD}{{\mathcal D}}
\newcommand{\cE}{{\mathcal E}}
\newcommand{\cM}{{\mathcal M}}
\newcommand{\cA}{{\mathcal A}}

\newcommand{\lip}{{\mathrm{Lip}^1}}
\newcommand{\lipl}{\mathrm{Lip}^1_{\mathrm{loc}}}

\newcommand{\One}{\mathds{1}}

\begin{document}

\title{A dual characterization of length spaces with application to Dirichlet metric spaces}
\author{Peter Stollmann}
\date{07-03-2009}
\maketitle

\begin{abstract}
We show that under minimal assumptions, the intrinsic metric induced by a strongly local Dirichlet form induces a length space. A main input is a dual characterization of length spaces in terms of the property that the 1-Lipschitz functions form a sheaf.
\end{abstract}

\section*{Introduction}
A metric space $(X,d)$ is said to be a \textbf{length space} or \textbf{path metric space}, whenever any two points in this space can be joined by a path with length arbitrarily close to the distance of these points. One main result of the present paper is that the intrinsic metric coming from a Dirichlet form gives a length space, whenever this intrinsic metric makes proper sense. This should be seen as complementary to results by Sturm \cite{Sturm-95b}, who considered the intrinsic metric under the additional assumption that all closed balls are compact and proved that then one even gets minimizing geodesics, i.e., paths with length equal to the distance. The classical energy form on open subsets of euclidean space shows that in general one will not encounter this more restrictive condition. 

In the investigation of the intrinsic metric of Dirichlet forms it turned out that a certain dual object of $X$, $\lip$, the space of 1-Lipschitz functions plays a central role: More precisely, the question whether $\lip$ is a \textbf{sheaf} will be crucial for the path metric property of a metric space. Here, we say that $\lip$ is a sheaf, if every function that is locally 1-Lipschitz is already globally 1-Lipschitz.

After introducing the necessary notions  in the following section, we turn to open subsets of euclidean space as a class of examples for which we can already illustrate the main questions and ideas. Then we prove the above mentioned dual characterization of length spaces and the final section is devoted to the proof of the fact that Dirichlet metric spaces are length spaces.

\textbf{Acknowledgment:} It is a pleasure to thank Ivan Veseli\'c who introduced me to the question that is answered here. Moreover, fruitful discussions with Maria are gratefully acknowledged. 

\section{Basic notions}
A metric $d$ on a set $X$ is a mapping $d:X\times X\to [0,\infty)$ such that $d(x,y)=0$ if and only $x=y$, that satisfies the triangle inequality $d(x,z)\le d(x,y)+d(y,z)$. If $d$ is allowed to take values in $[0,\infty]$ we speak of a \textbf{metric in the wide sense}. We write 
$$ B(x,r):=\{ y\in X \mid d(x,y)\le r\}\mbox{  and  }
U(x,r):=\{ y\in X \mid d(x,y) < r\}
$$
for the \textbf{closed} and \textbf{open balls}, respectively. A continuous mapping $\gamma:I\to X$ from an interval $I\subset \RR$ to $X$ is called a \textbf{path}. A metric in the wide sense induces a \textbf{length structure} in terms of
$$
L(\gamma):=\sup\{ \sum_{k=0}^{N}d(\gamma(t_{k+1}),\gamma(t_{k}))\mid t_0<t_1<\ldots < t_n, t_0,\ldots ,t_n\in I\} .
$$ 
The \textbf{path metric} induced by $d$ is given by 
$$
d_\ell(x,y):=\inf\{ L(\gamma)\mid \gamma : I\to X\mbox{  a path  }x,y\in \gamma(I)\}\in [0,\infty] ,
$$
with the usual convention that $\inf\emptyset =\infty$. 
See \cite{Gromov-07} for an axiomatic treatment of length structures that are not based on an a priori given metric.
The triangle inequality gives that $d(x,y)\le d_\ell(x,y)$ and we say that $(X,d)$ is a \textbf{length space}, provided that $d=d_\ell$. 
 
A canonical dual object for metric space is the space of 1-Lipschitz functions. Here are the respective notions: For $U\subset X$ denote
$$
 \lip (U,d) := \lip (U) := \{ f:X\to\RR |\; | f(x)-f(y)| \le d(x,y) \mbox{  for all  }
x,y\in U\} ,
$$
$$\lipl(U,d) := \lipl(U) := \{ f:X\to\RR \mid\forall x\in X\exists V\mbox{ open },
$$
$$
x\in V\mbox{  such that  }f|_{V\cap U}\in \lip (V\cap U) \} .
$$
Note that $\lip(U)\subset \lipl(U)$ but both are different in general, as can easily seen by considering a non-connected set $U$. We  say that $\lip$ \textbf{is a sheaf on} $U$, if 
$\lip(U) = \lipl(U)$. If that holds for $X=U$ we say that $\lip$ \textbf{is a sheaf}. Our first main result, Theorem \ref{sheaf} says that a locally complete space is a length space if and only if  $\lip$ is a sheaf.
It is clear that $d$ can be written as
$$
d(x,y) =\sup\{ f(y)-f(x)\mid f\in \lip (V)\}
$$
for every subset $V\subset X$ and all $x,y\in V$. The meaning of this observation will become clear immediately, when we discuss the intrinsic metric. 

\subsection*{Dirichlet metric spaces.}
The main application of our dual characterization and the starting point for the present paper are metrics induced by a strongly local, regular Dirichlet form, a notion we now briefly introduce. The starting setup is a locally compact metric space $X$ endowed with a regular Borel measure $m$ and a Dirichlet form $\cE$ in $L^2(X,m)$. We refer to \cite{BouleauH-91,Fukushima-80,FukushimaOT-94,MaR-92} for a thorough treatment of Dirichlet forms, a notion that goes back to \cite{BeurlingD-58, BeurlingD-59}. The example one should keep in mind is the classical Dirichlet form
$$
\cD(\cE):=W^{1,2}_0(\Omega), \cE[u,v]=\int_\Omega\nabla u(x)\cdot \nabla v(x) dx ,
$$
where $\Omega$ is an open subset of $\RR^d$ and $dx$ denotes integration with respect to Lebesgue measure. A Dirichlet form is called \textbf{regular}, if its domain $\cD(\cE)$ is dense both in $(C_c(X),\|\cdot\|_\infty)$ and in $(\cD(\cE),\|\cdot\|_\cE)$, where the \textbf{energy norm}  $\|\cdot\|_\cE$ is defined by $\| u \|^2_\cE:=\cE[u,u] +\| u\|^2$ and can be thought of as an analogue of the first order Sobolev norm, which appears for the classical Dirichlet form, as the reader will immediately notice.

If a regular Dirichlet form is \textbf{strongly local}, i.e., if
$$
 \cE[u,v]=0 \mbox{  for  }u\mbox{  constant on  }\supp v ,
$$
then $\cE$ can be represented in a way quite similar to the classical Dirichlet form. Namely, there exists a bilinear mapping $\Gamma$ from $ \cD(\cE)$ to the set of signed Radon measures $\cM_R$ such that
$$
\cE[u,v]=\int_Xd\Gamma(u,v) .
$$
This so called \textbf{energy measure} or \textbf{Lagrangian} can be defined via
$$
\int_Xd\Gamma(u,v) = \cE[u,\phi u]-\frac12\cE[u^2,\phi ]\mbox{  for  }u,\phi\in 
\cD(\cE)\cap C_c(X) .
$$ 
It can be extended to 
$$
\cD_{loc}:=\{ f\in L^2_{loc}(X)\mid \forall K\subset X\mbox{ compact }\exists u\in\cD(\cE):
u|_K =f|_K\}
$$
and inherits several important properties of the underlying Dirichlet form. E.g., $\Gamma$ is strongly local as well, meaning that, for open $V\subset X$ and $f\in \cD_{loc}$:
$$
\One_Vd\Gamma(f,f)=0\mbox{  whenever  }f\mbox{  is constant on  }V .
$$
Given the energy measure $\Gamma$, we can finally define the \textbf{intrinsic metric} in the following way: Consider
$$
\cA^1:=\{ f\in \cD_{loc}\cap C(X)\mid \Gamma(f,f)\le m\} ,
$$
where the latter inequality signifies that $ d\Gamma(f,f)$ is absolutely continuous with respect to the underlying measure $m$ with Radon-Nikod\'ym boundary bounded by $1$. In analogy with the classical Dirichlet form, $\cA^1$ can be thought of as those continuous functions, for which the gradient is bounded by $1$ in norm. We set
$$
d_\Gamma(x,y):= \sup\{ f(y) - f(x)\mid f\in \cA^1\}\in [0,\infty]
$$
and call it the \textbf{intrinsic metric} induced by $\cE$, see \cite{BiroliM-95, Sturm-94b,Sturm-95b}; some properties of the set $\cA^1$ can be found in the appendix of \cite{BoutetLS-08}. We will always assume that $\cE$ is \textbf{strictly local}, by which we mean that $d_\Gamma$ above is a metric in the wide sense and induces the original topology on $X$. Note that $d_\Gamma(x,y)=\infty$ occurs naturally if $x$ and $y$ are in different connected components of $X$, as was also discussed in \cite{Sturm-95b}. 

\section{The classical Dirichlet form on open subsets of $\RR^d$.}
Throughout this section $\Omega$ denotes an open subset of $\RR^d$. We will consider the usual euclidean metric $\rho$ as well as the intrinsic metric $d_\Gamma$ induced by the classical Dirichlet form defined above. In that case
$$
\Gamma(f,f)=| \nabla f(x)|^2dx \mbox{ for } f\in W^{1,2}_{loc}(\Omega) ,
$$
and 
$$
\cA^1=\{ f\in W^{1,2}_{loc}(\Omega)\cap C(\Omega) \mid | \nabla f(x)|\le 1 \; a.e. \} .
$$

\begin{remark}
 $(\Omega,\rho)$ is complete if and only if $\Omega=\RR^d$.
\end{remark}
Clearly, $\Omega=\RR^d\setminus \{ 0\}$ gives an example of a length space that is not geodesic, as open subsets are geodesic with respect to the euclidean metric if and only if they are convex. 
\begin{example}
 Let $\Omega=\RR^2\setminus(\{ 0\}\times [-1,1])$. Then $(\Omega,\rho)$ is not a length space. It is also easy to see that $\lip(\Omega)\varsubsetneq \lipl(\Omega)$. Choose, e.g.,
$f(x,y):=(1-|(x,y)|)_+$ for $(x,y)$ in the left halfplane and $f(x,y):=0$ on the right halfplane.
\end{example}
Of course, this example doesn't come as a surprise in view of Theorem \ref{sheaf} below. We now relate the intrinsic metric $d_\Gamma$ to the euclidean metric and obtain an explicit formula. We are not aware of any reference for this simple fact:
\begin{proposition}
 For any open subset $\Omega$ of $\RR^d$ we have that $d_\Gamma=\rho_\ell$, the length metric coming from euclidean distance.
\end{proposition}
\begin{proof}
 Pick $x,y\in \Omega$. If $x$ and $y$ lie in different connected components $U$ and $V$ of $\Omega$ then $\rho_\ell(x,y)=\infty$, since there is no path joining the two points. But $d_\Gamma(x,y)=\infty$ as well, as can be seen from picking $f_n:=n\One_V\in \cA^1$ as trial functions:
$$
 d_\Gamma(x,y)\ge \sup_{n\in\NN}(f_n(y)-f_n(x))=\infty .
$$
If $\rho_\ell(x,y)<\infty$, we can find $r>0$ such that $B(x,r)\subset \Omega$. Since the ball is convex we see that 
$$
| \rho_\ell(x,y_0)- \rho_\ell(x,y_1)|\le |y_0-y_1|\mbox{  for  }  
y_0,y_1\in B(x,r) ,
$$
which implies that $f:\Omega\to\RR,f(y):=\rho_\ell(x,y)$ for $ \rho_\ell(x,y)<\infty$ and $0$ else, defines a function in $\cA^1$. This gives
$$
d_\Gamma(x,y)\ge f(y)-f(x)=\rho_\ell(x,y) .
$$
Conversely, let $f\in \cA^1$ and $\gamma:[0,1]\to \Omega$ be a polygonal path from $x$ to $y$. Then

\begin{eqnarray*}
f(y)-f(x)&=&f(\gamma(1))-f(\gamma(0))\\
&=&\int_0^1f'(\gamma(t))\gamma'(t)dt \\
&\le&  \int_0^1 \gamma'(t)dt \\
&=& L(\gamma)
\end{eqnarray*}
By taking the $\inf$ over all path $\gamma$, we see that
$$
f(y)-f(x)\le \rho_\ell(x,y) .
$$
This gives $d_\Gamma(x,y)\le \rho_\ell(x,y)$.
\end{proof}
Since we have now calculated $d_\Gamma$ explicitly, we can record some simple consequences:
\begin{remark}
 Like for the euclidean distance, $(\Omega,d_\Gamma)$ is complete if and only if $\Omega=\RR^d$. For arbitrary $\Omega$, $(\Omega,d_\Gamma)$ is a length space, since $(\rho_\ell)_\ell=\rho_\ell$ for any metric; see \cite{Gromov-07}, 1.6.
\end{remark}
This gives a family of natural examples for which Dirichlet metrics are locally complete but not complete. In particular, our result Theorem \ref{dir} below, that all Dirichlet metrics define length spaces cannot be obtained from some Hopf-Rinow type result, as was done in \cite{Sturm-95b} under the additional assumption of completeness. See also \cite{Gromov-07,Rinow-61}. 
\section{Path metric property and $\lip$.}
In order to find a path with length close to the distance between $x$ and $y$ in a metric space one has to manage to find \textbf{approximate midpoints}, i.e., for every $x$ and $y$ and $\varepsilon>0$ there should be a $z$ such that 
$$
d(x,z)\le \frac12 d(x,y)+\varepsilon \mbox{  and  }d(z,y)\le \frac12 d(x,y)+\varepsilon .
$$
We will see now that such a property follows from sheaf properties of $\lip$. In this section, $(X,d)$ denotes a metric space.
\begin{proposition}\label{midpoints}
 Let $U$ be an open subset of $X$ such that $\lip(U)=\lipl(U)$. Then, for $d(x,y)< r_1+r_2$ it follows that $B(x,r_1)\cap B(y,r_2)\not=\emptyset$. In particular, we find approximate midpoints in $U$. 
\end{proposition}
\begin{proof}
 We proceed by contraposition and assume that $B(x,r_1)\cap B(y,r_2)=\emptyset$. Let $\delta >0$, so that $B(x,r_1-\delta)\cap B(y,r_2-\delta)=\emptyset$. Consider the function
$$
f(z):= -((r_1-\delta)-d(x,z))_++
((r_2-\delta)-d(z,y))_+ .
$$
It is clear that $f(y)=r_2-\delta$, $f(x)=-(r_1-\delta)$. Moreover, $f\equiv 0$ on $W:= U\setminus(B(x,r_1-\delta)\cup B(y,r_2-\delta))$. Using the covering $U(x,r_1), U(y,r_2),W$ we see that $f\in\lipl(U)$. Since $\lip(U)=\lipl(U)$, we get   that
$d(x,y)\ge f(y)-f(x)= r_1+r_2 -2\delta$. As $\delta$ was arbitrary, this proves that 
$d(x,y)\ge 
r_1+r_2$, as we wanted to show.
\end{proof}

Although the proof is clearly inspired by the proof of Lemma 3 in \cite{Sturm-95b}, both the assumption and the assertion are in fact quite different. It could also and will be used for Dirichlet metric spaces, where the defining function class $\cA^1$ is a sheaf. Here is our first main result:

\begin{theorem}\label{sheaf}
If $(X,d)$ is a length space then $\lip$ is a sheaf. Conversely, if $(X,d)$ is  locally complete and $\lip$ is a sheaf then $(X,d)$ is a length space.
\end{theorem}
\begin{proof}
 Assume that $(X,d)$ is a length space and let $f\in\lipl$. Let $x,y\in X$ and $\varepsilon>0$. Then there is a path $\gamma :[0,1]\to X$ such that $\gamma(0)=x$ and $\gamma(1)=y$ and $L(\gamma)\le d(x,y)+\varepsilon$. For every $t\in [0,1]$ there is an open (in $[0,1]$) interval $I_t$ containing $t$ and an open neighborhood $U_{\gamma(t)}$ such that 
$$\gamma(I_t)\subset U_{\gamma(t)}\mbox{ and  }f|_{U_{\gamma(t)}}\in \lip(U_{\gamma(t)}) .
$$ Since $[0,1]$ is compact, there are finitely many $0=t_1<\ldots <t_m=1$ such that 
$$
[0,1]\subset V_{t_1}\cup \ldots \cup V_{t_m}.
$$
and
$$
\gamma([0,1])\subset U_{\gamma(t_1)}\cup \ldots \cup U_{\gamma(t_m)}.
$$
Pick $s_1, \ldots s_{m-1}\in [0,1]$ such that 
$$
0\le t_1< s_1 <t_2<\ldots s_{m-1}<t_m=1\mbox{  with  } s_j\in I_{t_j}\cap I_{t_{j+1}} .
\mbox{  for  }j=1,\ldots ,m-1 .
$$
Then 
$$\gamma(s_j)\in  U_{\gamma(t_j)}\cap U_{\gamma(t_{j+1})} .
$$
Using the triangle inequality and the 1-Lipschitz property in the appropriate neighborhoods, we get:
\begin{eqnarray*}
 |f(x) - f(y)| &\le & \sum_{j=1}^{m-1}\left(|f(\gamma(t_j))-f(\gamma(s_j))|+
|f(\gamma(s_j))-f(\gamma(t_{j+1}))|\right) \\
&\le& \sum_{j=1}^{m-1}\left(d(\gamma(t_j),\gamma(s_j))+d(\gamma(s_j),\gamma(t_{j+1}))\right) \\
&\le& L(\gamma)\\
&\le& d(x,y) +\varepsilon .
\end{eqnarray*}
Since $\varepsilon>0$ was arbitrary, $|f(x)-f(y)|\le d(x,y)$.

To prove the converse, we want to show that, for fixed  $x\in X$, the function
$$
f:X\to\RR, f(y)=d_\ell(x,y)
$$
belongs to $\lipl$. If this is accomplished, we get that $f\in\lip$ using our assumption that $\lip$ is a sheaf, and  $f\in\lip$ gives the desired inequality:
$$
d(x,y)\ge |f(y)-f(x)|=d_\ell(x,y) .
$$
The first thing we have to check is that $f$ is properly defined, namely that $d_\ell(x,y)<\infty$ for all $y\in X$. The crucial step in this direction will also establish the local 1-Lipschitz property; denote $X_0:= \{ y\in X\mid d_\ell(x,y)<\infty\}$:

\noindent{\textbf{Claim:}} If $d_\ell(x,y)<\infty$, there exists $r>0$ such that 
$B(y,r)\subset X_0$ and for all
$y_0,y_1\in B(y,r)$:
$$
|d_\ell(x,y_0)-d_\ell(x,y_1)|\le d(y_0,y_1) . 
$$  
We have to find $r>0$ such that every  two points $y_0,y_1\in B(y,r)$ can be joined by a path with length arbitrarily close to $d(y_0,y_1)$, since this obviously implies the claim.
The existence of approximate midpoints settled in the preceding Proposition plus completeness (which we do have locally) will allow us to construct the desired path. This conclusion is well established, see \cite{Gromov-07}, Theorem 1.8. 
So take $r>0$ so small that $B(y,42r)$ is complete and $0<\varepsilon_1 <\frac12 r\wedge \frac12\varepsilon$; let $y_0,y_1\in B(y,r)$ and $\delta:=d(y_0,y_1)$.
From Proposition \ref{midpoints} we get a point $y_\frac12\in B(y_0,\frac12 \rho +\frac12\varepsilon_1)\cap B(y_1,\frac12 \rho +\frac12\varepsilon_1)$; in particular, $y_\frac12\in B(y,(1+\frac34)r)$; proceeding by induction, choosing a sequence of
$\varepsilon_k$ that decays rapidely enough, like in the above mentioned reference, we find a map 
$$
\gamma:\left\{ \frac{k}{2^n}\mid n\in\NN, k\in\NN_0, k\le 2^n\right\}\to B(y,42r)
$$
that is uniformly continuous with
$$
d(\gamma(\frac{k}{2^n}), \gamma(\frac{k+1}{2^n}))\le \delta\frac{1}{2^n}(1+\varepsilon) ,
$$
for all $k\le 2^n-1$. This implies that $d(\gamma(p),\gamma(q))\le \delta(1+\varepsilon)|p-q|$ for each pair of dyadic rationals in $[0,1]$. By completeness of the ball $B(y,42r)$, $\gamma$ extends to a $(1+\varepsilon)$-Lipschitz continuous path with length bounded by  $\delta(1+\varepsilon)$. 

From the claim we now get that $X_0$ is open and closed. Since $x\in X_0$, $X_0$ is nonempty and so must agree with $X$. In fact, as we already observed above, $\lip$ can only be a sheaf if the underlying space is connected. The claim also gives $f\in\lipl$, which completes the proof.
\end{proof}

\section{Strictly local Dirichlet spaces are length spaces.}

We now consider the setup introduced in Section 2 above: $X$ is a locally compact space,
$\cE$ a strictly local Dirichlet form and so comes with an energy measure $\Gamma$, for which
$$
\cA^1=\{ f\in \cD_{loc}\cap C(X)\mid \Gamma(f,f)\le m\} 
$$
separates the points of $X$ and
$$
d(x,y)=d_\Gamma(x,y)= \sup\{ f(y) - f(x)\mid f\in \cA^1\}\in [0,\infty]
$$
defines a metric in the wide sense that induces the original topology on $X$. In particular, small enough balls will be compact and hence complete. Clearly,
$$
\cA^1\subset \lip .
$$
Since also
$$
\lipl \subset \cA^1 ,
$$
as was established in \cite{LenzW-pre}, we almost have the following result:

\begin{theorem}
\label{dir} If $d$ is the intrinsic metric of a strictly local Dirichlet form $\cE$ on $X$, then $(X,d)$ is a length space.  
\end{theorem}
\begin{proof}
 The only difficulty we have to overcome is the fact that both $d$ and the corresponding path metric $d_l$ may take the value $\infty$. But, as we see, that happens simultaneously.
We fix $x\in X$.

\noindent{\textbf{1st Step:}} \emph{$d$ admits approximate midpoints}.

\noindent We can use the argument from the proof of Proposition \ref{midpoints}, keeping in mind that the property defining $\cA^1$ is local and so $\cA^1$ is a sheaf.

\noindent{\textbf{2nd Step:}} \emph{For every $y\in X$ there is $r>0$ such that
$$
d(y_0,y_1)=d_\ell(y_0,y_1)\mbox{  for all  }
y_0,y_1\in B(y,r) .
$$}
This follows from local completeness and the existence of  approximate endpoints exactly as in the proof of Theorem \ref{sheaf}. Again like in the latter proof we get that

\noindent{\textbf{3rd Step:}} \emph{The set $Y:=\{ y\in X\mid d_\ell(x,y)<\infty \}$ is open and closed.}

\noindent{\textbf{4th Step:}} \emph{If $d(x,y)<\infty$ then $d_\ell(x,y)<\infty$ as well. }
In fact, if $d_\ell(x,y)=\infty$, we have the $y\not\in Y$ and $x\in X_0$. Since the latter set is open and closed, it follows that $n\One_{Y}\in\cA^1$ (since $d\Gamma(\One_Y)=0$ by locality) for all $n\in\NN$. Therefore,
$$
d(x,y)\ge n\One_{Y}(x)-n\One_{Y}(y)=n
$$
for all $n\in\NN$.

\noindent{\textbf{5th Step:}} \emph{Let
$f:X\to\RR, f(y):= d_\ell(x,y)$ for $y\in Y$ and $0$ else. Then $f\in\cA^1$.}

By what we saw in the second step, $f\in\lipl$. Since  $\lipl \subset \cA^1$ as we mentioned above, we get the desired property of $f$.

We can now finish the proof as follows: The estimate $d(x,y)\le d_\ell(x,y)$ is clear. If $d_\ell(x,y)=\infty$, we know from the 4th Step that  $d(x,y)=\infty$ as  well. Thus it remains to consider $d(x,y)<\infty$, in other words $y\in Y$. But this gives
$$
d(x,y)\ge f(y)-f(x)=  d_\ell(x,y)$$ 
for $f$ as above.
\end{proof}
\begin{remark}
In the 1st step of the preceding proof we could also use a local version of Sturm's Lemma 3 from \cite{Sturm-95b}. Conversely, the property proved in the first step combined with compactness easily gives midpoints and this implies that one gets local minimizing geodesics.
\end{remark}

Apart from the examples coming from the classical Dirichlet form, there are many further classes of examples that fall into the framework covered by the preceding result. We mention second order elliptic \cite{BiroliM-95,MaR-92,Sturm-95a},
subelliptic operators \cite{FeffermanP-81,FeffermanSC-86,Jerison-86,JerisonSC-86,NagelSW-85} and  Laplace-Beltrami operators on manifolds. 

The above results were obtained in connection with a joint work on connectedness and irreducibility properties of Dirichlet forms with D. Lenz and I. Veseli\'c \cite{LenzSV-pre}.

{\sc Fakult\"at f\"ur Mathematik,
           Technische Universit\"at,\\ 09107 Chemnitz, Germany}\\
\texttt{{peter.stollmann@mathematik.tu-chemnitz.de }}
\end{document}